\documentclass[a4paper,12pt]{article}
\usepackage{tabularx}
\usepackage{array}
\usepackage{amsmath,amssymb}
\usepackage{enumerate}
\usepackage{epic}
\usepackage{bbm}
\usepackage{theorem}

\usepackage{color}
\usepackage[dvips,final]{graphicx} 

\newcommand{\NN}{\mathbb{N}}
\newcommand{\ZZ}{\mathbb{Z}}
\newcommand{\QQ}{\mathbb{Q}}
\newcommand{\RR}{\mathbb{R}}

\newcommand{\cF}{\ensuremath{\mathcal{F}}}

\newcommand{\<}{\ensuremath{\langle}}
\renewcommand{\>}{\ensuremath{\rangle}}

\newcommand{\supp}{\operatorname{supp}}
\newcommand{\ini}{\operatorname{in}}

\newtheorem{Theorem}{Theorem}[section]
\newtheorem{Lemma}[Theorem]{Lemma}
\newtheorem{Proposition}[Theorem]{Proposition}
\newtheorem{Corollary}[Theorem]{Corollary}

\newtheorem{Example}[Theorem]{Example}

\newenvironment{Proof}{\noindent{\it Proof.\/}}{\hfill $\square$\medskip}

\begin{document}
\title{The generic Gr\"obner walk}
\author{
K.~Fukuda\thanks{
Institute for Operations Research, ETH Z\"urich, Ch-8092, Switzerland,
fukuda@ifor.math.ethz.ch},\ 
A.~N.~Jensen\thanks{
Institut for Matematiske Fag,
Aarhus Universitet,
DK-8000 \AA rhus, Denmark, 
ajensen@imf.au.dk},\ 
N.~Lauritzen\thanks{
Institut for Matematiske Fag,
Aarhus Universitet,
DK-8000 \AA rhus, Denmark, 
niels@imf.au.dk},\ 
R.~Thomas\thanks{
Department of Mathematics, 
University of Washington,
Seattle, WA 98195, USA, 
thomas@math.washington.edu}
}

\maketitle

\begin{abstract}
  The Gr\"obner walk is an algorithm for conversion between Gr\"obner 
  bases for different term orders. It is based on the polyhedral 
  geometry of the Gr\"obner fan and involves tracking a line between 
  cones representing the initial and target term order.
  An important parameter is explicit numerical perturbation of this line. 
  This usually involves both time and space demanding arithmetic of 
  integers much larger than the input numbers.
  In this paper we show how the explicit line may be replaced by a formal line 
  using Robbiano's characterization of group orders on $\QQ^n$.
  This gives rise to the generic Gr\"obner walk involving only 
  Gr\"obner basis conversion over facets and computations with marked 
  polynomials. 
  The infinite precision integer arithmetic is replaced by term order 
  comparisons between (small) integral vectors. This makes it possible
  to compute with infinitesimal numbers and perturbations in a consistent 
  way without introducing unnecessary long integers.   
  The proposed technique is closely related to the lexicographic 
  (symbolic) perturbation method used in optimization and computational
  geometry. We report on computations with toric
  ideals, where a version of our algorithm in certain cases
  computes test sets for hard integer knapsack problems significantly faster 
  than the Buchberger algorithm.
  
\end{abstract}


\section{Introduction}

Let $R = k[x_1, \dots, x_n]$ denote the polynomial ring in $n$
variables over a field $k$. Gr\"obner basis computations in $R$ tend
to be very expensive for certain term orders (like the lexicographic
order). Therefore it often pays to compute Gr\"obner bases for
``easier'' term orders and convert them into Gr\"obner bases for the
desired term order. For zero-dimensional ideals this can be
accomplished by the FGLM-algorithm \cite{FGLM}. For general ideals the
Gr\"obner walk algorithm \cite{CKM} can be applied. 

Let $\prec_1$ and $\prec_2$ be term orders on $R$.  The usual
Gr\"obner walk proceeds from the reduced Gr\"obner basis $G$ for $I$
over $\prec_1$ by tracking a line $\omega(t) = (1-t) \omega_0 + t
\tau_0, \,\, 0\leq t \leq 1$, where $\omega_0$ and $\tau_0$ are
vectors in the respective Gr\"obner cones $ C_{\prec_1}(I)$ and
$C_{\prec_2}(I)$ of $I$. At $t = 0$ the Gr\"obner basis is known.  The
line $\omega(t)$ is tracked through the Gr\"obner fan of $I$ and
Gr\"obner bases are computed at common faces of successive Gr\"obner
cones. At $t = 1$ we reach the reduced Gr\"obner basis for $I$ over
$\prec_2$.

The efficiency of the Gr\"obner walk rests on clever choices
of $\omega_0$ and $\tau_0$. A choice of $\omega_0$ and $\tau_0$ on low
dimensional faces of Gr\"obner cones may lead to very heavy Gr\"obner
basis calculations along $\omega(t)$. Often (but not always) it pays
to choose $\omega_0$ and $\tau_0$ generically inside $C_{\prec_1}(I)$
and $C_{\prec_2}(I)$ and ensure that $\omega(t)$ only intersects
common faces of low codimension on its way to the target term order
$\prec_2$.

The initial reduced Gr\"obner basis $G$ over $\prec_1$ makes it
possible to compute an interior point in $C_{\prec_1}(I)$. Computing
an interior point in the target cone $C_{\prec_2}(I)$ is considerably
more difficult, since we do not know the reduced Gr\"obner basis over
$\prec_2$ in advance.  Tran \cite{Tran} approached this problem using
general degree bounds on polynomials in Gr\"obner bases.  The general
degree bounds in Tran's approach may lead to integral weight vectors
with $10,000$-digit entries in representing a lexicographic interior
point in the case of polynomials of degree $10$ in $10$ variables.

In this paper we give an algorithm where the line $\omega(t)$ is
replaced by a (formal) line $\Omega(t)$ between suitably chosen
perturbations given by $\prec_1$ and $\prec_2$ and $I$. It turns 
out that the numerical dependence on $I$ disappears in our algorithm
and that $\Omega(t)$ may 
be viewed as a line which can be used for all ideals in the Gr\"obner
walk from $\prec_1$ to $\prec_2$. The formal line has the property
that its initial and target points are always in the interior of
the Gr\"obner cones. Furthermore the common faces that $\Omega(t)$
intersect are all facets.

In the classical Buchberger algorithm \cite{Buch} for computing
Gr\"obner bases one only computes with term orders and initial
terms of polynomials.
Tracking $\Omega(t)$ gives a ``Buchberger-like'' Gr\"obner walk
algorithm, where one only needs to compute with marked polynomials and
term orders.  On smaller examples the algorithm can easily be carried
out by hand (cf.~\S \ref{Section:Example}).  

We have observed
some interesting experimental results using a
version of the generic walk tailored to lattice ideals \cite{L}. 
When the generic walk is applied in computing full test sets 
for feasibility of the hard integer knapsacks from \cite{AL}, the
natural initial and target vectors are rather close in the
Gr\"obner fan. This leads to very fast computations of
test sets. These examples 
with polynomials of high degree in many variables seem
out of reach for the classical Gr\"obner walk.
We report on computational experiments
in the last section of this paper.

An understanding of our algorithm requires a firm grip on the usual
Gr\"obner walk algorithm. Therefore \S \ref{preliminaries} and \S
\ref{Section:GW} recalls and proves fundamental results for the usual
Gr\"obner walk using which we transition to the generic Gr\"obner walk
in \S \ref{generic walk}.

The basic technique we propose to avoid explicit perturbation
is not quite new.  The key idea of implicit (symbolic) perturbation was 
proposed by Charnes in 1952 to make Dantzig's simplex method 
for linear programming finite.   The method is now known as
the lexicographic perturbation method, see \cite[page 34]{c-lp-83},
and used by many reliable implementations of the simplex method.
In computational geometry, similar symbolic perturbation schemes
are used to treat input data points in $\RR^n$ as if they were in
general position, see \cite[page 14]{bkos-cgaa-00}.

\section{Preliminaries} \label{preliminaries}

In this section we recall the basics of convex polyhedral cones. We
emphasize a crucial result from the theory of group orders (Lemma
\ref{Lemma:Interior}) and recall the construction of the (restricted)
Gr\"obner fan.

\subsection{Cones and fans}
\

A {\it convex polyhedral cone\/} is a set
$$
C(v_1, \dots, v_r) = \RR_{\geq 0} v_1 + \cdots + \RR_{\geq 0} v_r\subseteq
\RR^n
$$
where $v_1, \dots, v_r\in \RR^n$. In the following a {\it cone\/} 
will refer to a convex polyhedral cone. The dual of a 
cone $C\subseteq \RR^n$ is  
$$
C^\vee = \{\omega\in \RR^n \mid \<\omega, v\> \geq 0, \mathrm{\ for\ 
  every\ }v\in C\}.
$$
The dual of a cone is a cone and the intersection of two cones is a
cone.  The {\it dimension\/} of a cone is the dimension of the linear
subspace it spans.  For a vector $u\in\RR^n$ we let $u^\perp = \{x\in
\RR^n\mid \<u, x\> = 0\}$. A face $F \subseteq C$ of a cone $C$ is a
subset $F = u^\perp\cap C$, where $u\in C^\vee$.  Faces of codimension
one in $C$ are called {\it facets\/}.

A collection $\cF$ of cones and their faces is called a {\it fan\/} if
for every $C_1, C_2\in \cF$ we have $C_1\cap C_2\in \cF$ and $C_1\cap
C_2$ is a common face of $C_1$ and $C_2$.

\subsection{Rational group orders on $\QQ^n$}\label{RTO}
\ 

Let $(A, +)$ be an abelian group. Recall that a {\it group order\/}
$\prec$ on $A$ is a total order $\prec$ on $A$ such that
$$
x \prec y \implies x + z \prec y + z
$$
for every $x, y, z\in A$.

Let $\omega = (\omega_1, \dots, \omega_n)\subset \QQ^n$ be a
$\QQ$-vector space basis for $\QQ^n$.  Then we get a group order
$\prec_\omega$ on $\QQ^n$ given by $u \prec_\omega v$ if and only if
$$
(\<\omega_1, u\>, \dots, \<\omega_n, u\>) <_{\mathrm{lex}} 
(\<\omega_1, v\>, \dots, \<\omega_n, v\>),
$$
where $<_{\mathrm{lex}}$ refers to the lexicographic order on
$\QQ^n$. We call such a group order {\it rational\/}. To describe arbitrary
group orders on $\QQ^n$ similarly, one needs a more general setup
including real vectors (see \cite{Robbiano}). To ease the exposition
we will restrict ourselves to rational group orders. A group order
refers to a rational group order in the following.  For a rational
$\epsilon > 0$ we put
$$
\omega_\epsilon = \omega_1 + \epsilon \omega_2 + \cdots +
\epsilon^{n-1} \omega_n.
$$
The following well known lemma plays a key role in the 
generic Gr\"obner walk.

\begin{Lemma}\label{Lemma:Interior}
  Let $\omega = (\omega_1, \dots, \omega_n)\subset \QQ^n$ be a
  $\QQ$-basis.  Suppose that $F\subset \QQ^n$ is a finite set of
  non-zero vectors with $0\prec_\omega v$ for $v\in F$. Then there
  exists $0 < \delta\in \QQ$ such that $\<\omega_\epsilon, v\> > 0$
  for every $v\in F$ and $\epsilon\in \QQ$ with $0 < \epsilon <
  \delta$.
\end{Lemma} 
\begin{Proof}
  We prove this by induction on $n$. The case $n = 1$ is clear.  For
  $n > 1$ we may find $0 <\delta_0\in \QQ$ such that
$$
\<\omega_{n-1} + \epsilon \omega_n, v\> > 0
$$
for every $v\in F$ with $\<\omega_{n-1}, v\> > 0$ and $\epsilon\in
\QQ$ with $0 < \epsilon < \delta_0$. Therefore $0\prec_{\omega'} v$
for $\omega' = (\omega_1, \dots, \omega_{n-2}, \omega_{n-1} + \epsilon
\omega_n)$ for every $v\in F$ if $0 < \epsilon <\delta_0$.

By induction there exists $0 < \delta_1\in \QQ$ such that
$\<\omega'_\epsilon, v\> > 0$ for every $v\in F$ and $\epsilon\in \QQ$
with $0 < \epsilon < \delta_1$. Putting 
$\delta = \min(\delta_0, \delta_1)$ we get $\<\omega_\epsilon, v\> >
0$ for every $v\in F$ and $\epsilon\in \QQ$ with $0 < \epsilon < \delta$.
\end{Proof}

A group order $\prec$ on $\QQ^n$ is called a term order if $0\prec v$
for every $v\in \NN^n$. This is equivalent to $0 \prec e_i$ where
$e_i$ denotes the $i$-th canonical basis vector for $i = 1, \dots, n$.
As a consequence of Lemma \ref{Lemma:Interior} we get the following
corollary.

\begin{Corollary}\label{Corollary:FullDim}
  Let $F\subset \QQ^n$ be a finite set of positive vectors for the
  group order $\prec$ i.e. $v \succ 0$ for every $v\in F$. Then there
  exists $\omega\in \QQ^n$ such that
$$
\<\omega, v\> > 0
$$
for every $v\in F$. If $\prec$ is a term order, we may assume that
$\omega$ has positive coordinates.
\end{Corollary}

\subsection{The Gr\"obner fan}\label{GF}
\

Let $R = k[x_1, \dots, x_n]$ denote the ring of polynomials in $n$
variables over a field $k$. It is convenient to view $R$ as the
semigroup ring $k[\NN^n]$. We briefly recall the construction of the
(restricted) Gr\"obner fan (cf.~\cite{MR}) for an arbitrary ideal in
$R$.

Fix a group order $\prec$ on $\QQ^n$.  For a polynomial $f =
\sum_{v\in \NN^n} a_v x^v\in R$ we let $\supp(f) = \{v\in \NN^n \mid
a_v \neq 0\}$ and $\ini_\prec(f) = a_u x^u$, where $u =
\max_\prec\supp(f)$. For a vector $\omega\in \RR^n$ we let
$\ini_\omega(f)$ denote the sum of terms $a_v x^v$ in $f$ maximizing
the {\it $\omega$-weight\/} $\<\omega, v\>$.  We call $f$
$\omega$-{\it homogeneous\/} if $f = \ini_\omega(f)$.  To an ideal
$I\subseteq R$ we associate the ideals $\ini_\prec(I) =
\<\ini_\prec(f)\mid f\in I\setminus\{0\}\>$ and $\ini_\omega(I) =
\<\ini_\omega(f)\mid f\in I\>$. These ideals may be viewed as
deformations of the original ideal $I$. The {\it initial ideal\/}
$\ini_\prec(I)$ is generated by monomials. This does not hold for
$\ini_\omega(I)$ in general (unless $\omega$ is chosen generically).
  
Now define
$$
\partial_\prec(f) = \{u - u' \mid u'\in \supp(f)\setminus\{u\}\} \subset
\ZZ^n,
$$
where $a_u x^u = \ini_\prec(f)$. For a finite set $F\subseteq R$ of 
polynomials we let
$$
\partial_\prec(F) = \bigcup_{f\in F} \partial_\prec(f)
$$
and 
\begin{align*}
  C_\prec(F) &= C(\partial_\prec(F))^\vee \cap \RR^n_{\geq 0}\\
  & = \{\omega\in \RR^n_{\geq 0}\mid \<\omega, v\>\geq 0, v\in
  \partial_\prec(f), f\in F\}.
\end{align*}
Notice that $\dim C_\prec(F) = n$ by Corollary \ref{Corollary:FullDim}
and that
$$
C_\prec(F) = \{\omega\in \RR^n_{\geq 0}\mid
\ini_\prec(\ini_\omega(f)) = \ini_\prec(f)\mathrm{\ for\ every\ }f\in
F\}.
$$
A {\it Gr\"obner basis\/} for $I$
over $\prec$ is a finite set of polynomials $G = \{g_1, \dots, g_r\}
\subseteq I$ such that
$$
\<\ini_\prec(g_1), \dots, \ini_\prec(g_r)\> =
\ini_\prec(I).
$$
The Gr\"obner basis $G$ is called {\it minimal\/} if none of $g_1,
\dots, g_r$ can be left out and {\it reduced\/} if the coefficient of
$\ini_\prec(g_i)$ is $1$ and $\ini_\prec(g_i)$ does not divide any of
the terms in $g_j$ for $i\neq j$ and $i, j = 1, \dots, r$. A reduced
Gr\"obner basis is uniquely determined.  Minimal Gr\"obner bases exist
for arbitrary group orders.  However, Gr\"obner bases over arbitrary
group orders do not necessarily generate the ideal (as opposed to
Gr\"obner bases over term orders).  Similarly, the reduced Gr\"obner
basis is only guaranteed to exist for term orders.

To define the Gr\"obner fan we now specialize to the case where
$\prec$ is a term order.
The {\it Gr\"obner cone\/} $C_\prec(I)$ of an ideal $I$ over $\prec$ is
defined as $C_\prec(G)$, where $G$ is the reduced Gr\"obner basis of
$I$ over $\prec$. The {\it Gr\"obner fan\/} of $I$ is defined as the set of
cones $C_\prec(I)$ along with their faces, where $\prec$ runs through
all term orders. This is a finite collection of cones 
\cite[Theorem 1.2]{St} and one may
prove that it is a fan 
(Propositions 2.3 and 2.4 in \cite{St} give a proof assuming 
non-negative weight vectors).  The following proposition 
shows that $C_\prec(I)$ is the largest cone among $C_\prec(G)$, where 
$G$ is a Gr\"obner basis for $I$ over $\prec$.

\begin{Proposition}
Let $G$ be a (not necessarily reduced) Gr\"obner basis for $I$ 
over $\prec$. Then 
$$
C_\prec(G) \subseteq C_\prec(I).
$$
\end{Proposition}
\begin{Proof}
  Transforming $G$ into a minimal Gr\"obner basis $G'$ by omitting
  certain polynomials in $G$ clearly leads to an inclusion $C_\prec(G)
  \subseteq C_\prec(G')$. Transforming $G'$ into the reduced Gr\"obner
  basis proceeds by a sequence of reduction steps: suppose that $f_i,
  f_j\in G'$ and that a term $x^v$ in $f_j$ is divisible by
  $\ini_\prec(f_i)$.  Then $f_j$ is replaced by $f_j' = f_j -
  (x^v/\ini_\prec(f_i)) f_i$.  This reduction may introduce ``new''
  monomials which are not present in $f_j$. More precisely if $w\in
  \supp(f_j')$, then $w\in \supp(f_j)$ or $w = v - u + u'$, where $a_u
  x^u = \ini_\prec(f_i)$ and $u'\in \supp(f_i)$. In the latter case we
  get $w' - w = (w' - v) + (u - u')$, where $a_{w'} x^{w'} =
  \ini_\prec(f_j)$. Let $G''$ denote the Gr\"obner basis obtained by
  replacing $f_j$ with $f_j'$. Then $C(\partial_\prec(G')) \supseteq
  C(\partial_\prec(G''))$ and thereby $C_\prec(G')\subseteq
  C_\prec(G'')$. Since the reduced Gr\"obner basis is obtained using a
  finite number of these reduction steps we have proved the inclusion.
\end{Proof}

For a specific term order one may have infinitely many cones given by
different minimal Gr\"obner bases. As an example consider the ideal $I
= \<x, y\> \subset k[x, y]$. If $n$ is a positive natural number then
$G_n = \{x - y^n, y\}$ is a minimal Gr\"obner basis for $I$ over the
lexicographic order $\prec$ with $x \succ y$. In this case
$$
C_\prec(I) \supsetneq C_\prec(G_1) \supsetneq
C_\prec(G_2)\supsetneq \cdots.
$$

\section{The Gr\"obner walk}\label{Section:GW}

We outline the basic idea of the Gr\"obner walk \cite{CKM} and give a
new lifting step using reduction modulo the known Gr\"obner basis.

Let $\prec_1$ and $\prec_2$ be term orders and $I$ an ideal in $R$.
Suppose that we know the reduced Gr\"obner basis $G$ for $I$ over
$\prec_1$.  If
$$
\omega\in C_{\prec_1}(I) \cap C_{\prec_2}(I)
$$
lies on the common face of the two Gr\"obner cones, then $G_\omega
= \{\ini_\omega(g) \mid g\in G\}$ is the reduced Gr\"obner basis for
$\ini_\omega(I)$ over $\prec$.  Now a ``lifting'' of $G_\omega$ to a
Gr\"obner basis for $I$ over $\prec_2$ is required. The procedure for
this is based on Proposition \ref{Proposition:Walk} below. It involves
a Gr\"obner basis computation for $\ini_\omega(I)$ over $\prec_2$. The
point is that if $F=C_{\prec_1}(I)\cap C_{\prec_2}(I)$ is a high
dimensional face (like a facet) and $\omega$ is in the relative
interior of $F$, the ideal $\ini_\omega(I)$ is close to a monomial
ideal and this Gr\"obner basis computation becomes very easy.

Given a term order $\prec$ and a vector $\omega\in \RR_{\geq 0}^n$ we
define the new term order $\prec_\omega$ by $u \prec_\omega v$ if and
only if $\<u, \omega\> < \<v, \omega\>$ or $\<u, \omega\> = \<v,
\omega\>$ and $u \prec v$.  We record the following well known lemma.

\begin{Lemma} \cite[Proposition 1.8]{St} 
\label{Lemma:initomega}
Let $I\subseteq R$ be any ideal and $\omega\in \RR_{\geq 0}^n$. Then
$$
\ini_\prec(\ini_\omega(I)) = \ini_{\prec_\omega}(I).
$$
\end{Lemma}

The lifting step (Proposition \ref{Proposition:Walk}(ii) below) 
in the following proposition is different from the lifting step 
in the usual Gr\"obner walk \cite[Subroutine 3.7]{St}.

\begin{Proposition}\label{Proposition:Walk}
  Let $I\subseteq R$ be an ideal and $\prec_1, \prec_2$ term orders on
  $R$. Suppose that $G$ is the reduced Gr\"obner basis for $I$ over
  $\prec_1$.  If $\omega\in C_{\prec_1}(I)\cap C_{\prec_2}(I)$, then
\begin{enumerate}[(i)]
\item The reduced Gr\"obner basis for $\ini_\omega(I)$ over $\prec_1$
  is $G_\omega = \{\ini_\omega(g) \mid g\in G\}$.
\item If $H$ is the reduced Gr\"obner basis for $\ini_\omega(I)$ over
  $\prec_2$, then
$$
\{f - f^G\mid f\in H\}
$$
is a minimal Gr\"obner basis for $I$ over ${\prec_2}_\omega$. Here
$f^G$ is the unique remainder obtained by dividing $f$ modulo $G$.
\item The reduced Gr\"obner basis for $I$ over ${\prec_2}_\omega$
  coincides with the reduced Gr\"obner basis for $I$ over $\prec_2$.
\end{enumerate}
\end{Proposition}
\begin{Proof}
  Given a term order $\prec$ and a vector $\omega\in C_\prec(I)$, the
  reduced Gr\"obner bases for $I$ over $\prec$ and $\prec_\omega$
  agree. This proves (iii) and (i) taking Lemma \ref{Lemma:initomega}
  into consideration.  Suppose that $f$ is an $\omega$-homogeneous
  polynomial (cf.~\S \ref{GF}) in $\ini_\omega(I)$. Using the division
  algorithm in computing the unique remainder $f^G$, we keep reducing
  terms with the same $\omega$-weight as the terms in $f =
  \ini_\omega(f)$. Since $\ini_{\prec_1}(\ini_\omega(g)) =
  \ini_{\prec_1}(g)$ for $g\in G$ and $f^{G_\omega} = 0$, we see
  that all terms in $f^G$ will have $\omega$-weight strictly less than
  the terms in $f$. Therefore
$$
\ini_\omega(f) = \ini_\omega(f - f^G).
$$
Now suppose that $\{f_1, \dots, f_s\}$ is the reduced Gr\"obner
basis for $\ini_\omega(I)$ over $\prec_2$. In particular we get that
$f_i$ is $\omega$-homogeneous for $i = 1, \dots, s$. Then
\begin{align*}
{\ini_{\prec_2}}_\omega(I) &= \ini_{\prec_2}(\ini_\omega(I)) =
\<\ini_{\prec_2}(f_1), \dots, \ini_{\prec_2}(f_s)\>\\
&=\<\ini_{\prec_2}(\ini_\omega(f_1)), \dots,
\ini_{\prec_2}(\ini_\omega(f_s))\>\\
&=\<\ini_{\prec_2}(\ini_\omega(f_1 - f_1^G)), \dots,
\ini_{\prec_2}(\ini_\omega(f_s - f_s^G))\>\\
&=\<{\ini_{\prec_2}}_\omega(f_1 - f_1^G), \dots,
{\ini_{\prec_2}}_\omega(f_s - f_s^G)\>.
\end{align*}
This proves that $\{f_1 - f_1^G, \dots, f_s - f_s^G\}\subseteq I$ is a 
(minimal) Gr\"obner basis for $I$ over ${\prec_2}_\omega$.
\end{Proof}

Proposition \ref{Proposition:Walk} may be turned into a Gr\"obner
basis conversion algorithm as shown in the following section.

\subsection{Conversion along a line}

A natural approach to Gr\"obner basis conversion is to trace the line
between vectors in different Gr\"obner cones and update Gr\"obner
bases successively using Proposition \ref{Proposition:Walk}. This 
process is called the Gr\"obner walk
\cite{CKM}. A good reference for this procedure is \cite[\S4]{CLO},
which inspired the following. 
We sketch the first step of the Gr\"obner walk.  The succeeding
steps of the Gr\"obner walk are similar.  Suppose that $\omega_0\in
C_{\prec_1}(I), \tau_0\in C_{\prec_2}(I)$ and that $G$ is the reduced
Gr\"obner basis for $I$ over $\prec_1$.  Here $\prec_1$ and $\prec_2$
are rational term orders (cf.~\S\ref{RTO}) given by $\QQ$-bases
$\omega = (\omega_1, \dots, \omega_n)$ and $\tau = (\tau_1, \dots,
\tau_n)$ respectively.  Then we
consider the line
$$
\omega(t) = (1-t)\omega_0 + t \tau_0,\,\, 0\leq t \leq 1
$$
in the Gr\"obner fan of $I$ from $\omega_0$ to $\tau_0$.  Initially
we know the reduced Gr\"obner basis at $\omega(0) = \omega_0$ (being
G). Consider the ``last'' $\omega'=\omega(t')$ in $C_{\prec_1}(I) =
C_{\prec_1}(G)$. To be more precise $t'$ satisfies
\begin{enumerate}
\item
$0\leq t' < 1$
\item $\omega(t)\in C_{\prec_1}(I)$ for $t\in [0, t']$ and
  $\omega(t'+\epsilon)\not\in C_{\prec_1}(I)$ for every $\epsilon>0$.
\end{enumerate}
If no such $t'$ exists then $G$ is the reduced Gr\"obner basis over
$\prec_2$.  
If $t'$ exists $\omega(t')$ is on a proper face of $C_{\prec_1}(I)$
and $v\in \partial(G)$ exists with $\<\omega(t' + \epsilon), v\> < 0$
for $\epsilon > 0$.
This implies that $\<\tau_0,v\> < \<\omega_0, v\>$
and hence $\<\tau_0, v\> < 0$.  
 
This indicates the procedure for finding $t'$ given $G$. For $v\in
\partial(G)$ satisfying $\<\tau_0, v\> < 0$ we solve $\<\omega(t), v\>
= 0$ for $t$ giving
$$
t_v = \cfrac{\<\omega_0, v\>}{\<\omega_0, v\> - \<\tau_0, v\>}.
$$
Then $t'$ is the minimal among these $t_v$. In this case $\omega' =
\omega(t')$ lies on a proper
face $F$ of $C_{\prec_1}(I)$ and clearly
$$
\omega'\in C_{{\prec_2}_{\omega'}}(I).
$$
Now we use ${\prec_2}_{\omega'}$ as the term order $\prec_2$ in
Proposition \ref{Proposition:Walk}. The point is that we only need the
target term order $\prec_2$ to compute a Gr\"obner basis for
$\ini_{\omega'}(I)$ (not the notational beast ${\prec_2}_{\omega'}$).
The reason for this is that the Buchberger algorithm in this case
solely works with $\omega'$-homogeneous polynomials and ties are
broken with $\prec_2$.

To prove that we actually enter a new Gr\"obner cone we need to show
that $t' > 0$ (cf.~\cite{CLO}, \S5, (5.3) Lemma). In the initial step
it may happen that $t' = 0$. But if this is the case we may assume (in
the following step of the Gr\"obner walk) that $G$ is the reduced
Gr\"obner basis over $\prec = {\prec_2}_{\omega'}$. Since $\tau_0\in
C_{\prec_2}(I)$ is non-zero we may use $\tau_0$ as the first vector in
a $\QQ$-basis representing $\prec_2$. In this case assume that $t'=0$.
This means that we can find $v\in \partial(G)$ with $\<\omega', v\> =
0$ and $\<\tau_0, v\> < 0$ contradicting that $G$ is a Gr\"obner basis
over $\prec$.

We have outlined the procedure for tracking the line $\omega(t)$
through the Gr\"obner fan detecting when $\omega(t)$ leaves a cone.
The salient point of the generic Gr\"obner walk is that this
calculation can be done formally by choosing certain
generically perturbed $\omega_0$ and $\tau_0$ given by
$\prec_1$ and $\prec_2$.

Here are the steps of the usual Gr\"obner walk algorithm with the
modified lifting step. Recall that a marked polynomial is a polynomial
with a distinguished term, which is the initial term with respect to a
term order $\prec$. For a marked polynomial $f$, $\partial(f)$ is defined
in the natural way (cf.~the definition of $\partial_\prec(f)$ in \S \ref{GF}).
A marked Gr\"obner basis over a term order $\prec$ is a
Gr\"obner basis over $\prec$ with all initial terms (with respect to 
$\prec$) marked. For a marked Gr\"obner basis we let $\partial(G) = 
\cup_{f\in G} \partial(f)$.

\

\noindent
{\bf INPUT}: Marked reduced Gr\"obner basis $G$ for $I$ over a 
term order $\prec_1$, a term order $\prec_2$ along with 
$\omega_0\in C_{\prec_1}(I)$ and $\tau_0\in C_{\prec_2}(I)$.

\

\noindent
{\bf OUTPUT}: Reduced Gr\"obner basis for $I$ over $\prec_2$.

\begin{enumerate}[(i)]
\item
$t = -\infty$.
\item \label{Walk:Repeat}
{\bf Compute\_last\_t}. 
If $t = \infty$ output $G$ and halt.
\item
Compute generators $\ini_\omega(G) = \{\ini_\omega(g) \mid g\in G\}$ 
for $\ini_\omega(I)$ as
$$
\ini_\omega(g) = a^u x^u + \sum_{v\in S_g} a_v x^v,
$$ 
where $S_g = \{v\in \supp(g)\setminus \{u\} \mid t_{u - v} = t\}$ and
$a_u x^u$ is the marked term of $g\in G$.
\item
Compute reduced Gr\"obner basis $H$ for $\ini_\omega(I)$ over 
$\prec_2$ and mark $H$ according to $\prec_2$.
\item
Let
$$
H' = \{f - f^G \mid f\in H\}.
$$
Use marking of $H$ to mark $H'$.
\item
Autoreduce $H'$ and put $G = H'$.
\item
Repeat from (\ref{Walk:Repeat}).
\end{enumerate}

\noindent
{\bf Compute\_last\_t:}

\begin{enumerate}
\item
Let $V := \{v\in \partial(G) \mid \<\omega_0, v\> \geq  0
\mathrm{\ and\ }\<\tau_0, v\> < 0\mathrm{\ and\ } t \leq t_v\}$, where
$$
t_v = \cfrac{\<\omega_0, v\>}{\<\omega_0, v\> - \<\tau_0, v\>}.
$$
\item
If $V = \emptyset$, put $t = \infty$ and return. 
\item\label{Stepcrux} 
Let $t := \min\{t_v | v\in V\}$ and return.
\end{enumerate}

\section{The generic Gr\"obner walk} \label{generic walk}

In this section we show how certain generic choices of $\omega_0$ and
$\tau_0$ from \S \ref{Section:GW} lead to the
generic Gr\"obner walk algorithm. The crucial point is that step 
(\ref{Stepcrux}) of the procedure {\bf Compute\_last\_t} can be 
carried out formally using $\omega_0$ and $\tau_0$ from
well defined perturbations given the term orders $\prec_1$ and
$\prec_2$.

For an ideal $I\subseteq R$ we let $\partial(I)\subseteq \QQ^n$ denote
the union of $\partial_\prec(G)$, where $G$ runs through the finitely
many reduced Gr\"obner bases for $I$.  Let $\prec_1$ and $\prec_2$ be
two term orders given by $\QQ$-bases $\omega = (\omega_1, \dots,
\omega_n)$ and $\tau=(\tau_1, \dots, \tau_n)$ of $\QQ^n$ respectively.
Observe that $\omega_\eta$ and $\tau_\eta$ are in the interior of the
Gr\"obner cones $C_{\prec_1}(I)$ and $C_{\prec_2}(I)$ respectively for
sufficiently small positive $\eta$.  This follows from Lemma
\ref{Lemma:Interior}.  Now define
$$
C_{\prec_1, \prec_2} = \{v\in \RR^n \mid 0 \prec_1 v\mathrm{\ and\ 
} v \prec_2 0\}.
$$
Here $\prec_1, \prec_2$ are extended to group orders on $\RR^n$ using 
$\omega$ and $\tau$. 

\begin{Example}\label{Example2}
Suppose that $\prec_1$ is degree (reverse) lexicographic order and
$\prec_2$ lexicographic order with $y \prec_{1,2}x$. Then choosing 
$\omega = ((1, 1), (0, -1))$ and 
$\tau = ((1, 0), (0, 1))$,  we get  
$0 \prec_1 v$ imples $(0,0) <_{lex} (v_1+v_2, -v_2)$ and 
$v \prec_2 0$ implies $(v_1,v_2) <_{lex} (0,0)$. Intersecting the 
regions yielded gives (see Figure 1) 
$$
C_{\prec_1, \prec_2} = 
\{(x, y)\in \RR^2 \mid x+y > 0, x < 0\}. 
$$
\begin{figure}
\begin{center}
\includegraphics[height=90mm]{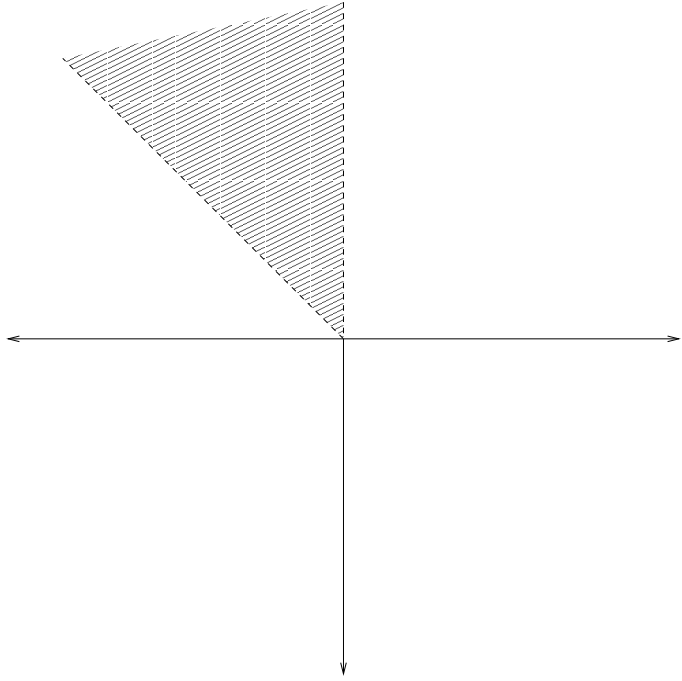}
\end{center}
\caption{$C_{\prec_1, \prec_2}$ for $\prec_1 =$degrevlex and $\prec_2$=lex}
\end{figure}
\end{Example}

To fully understand the choice of $\delta$ and $\epsilon$ in the
following we encourage the reader to compare with the computations in
($*$) and ($**$) below.  Define
$$
M_\tau = \{\<\tau_i, u\> v \mid i = 1, \dots, n; u, v\in \partial(I)\}.
$$
Corollary \ref{Corollary:FullDim} shows that there exists
sufficiently small positive $\delta$ such that
\begin{equation}\label{delta}
u \prec_1 v \iff \<\omega_\delta, u\> < \<\omega_\delta, v\>
\end{equation}
for $u, v\in M_\tau$. Suppose that $\delta$ satisfies (\ref{delta}).
Now put 
$$
N_\delta = \{\<\omega_\delta, u\> v \mid u, v\in \partial(I)\}.
$$
Again by Corollary \ref{Corollary:FullDim} we know that 
there exists sufficiently small positive $\epsilon$ such that
\begin{equation}\label{epsilon}
u \prec_2 v \iff \<\tau_\epsilon, u\> < \<\tau_\epsilon, v\>
\end{equation}
for $u, v\in N_\delta$. Suppose now that we pick $\delta$ according
to (\ref{delta}) and $\epsilon$ according to (\ref{epsilon}). If 
$v\in \partial(I)\cap C_{\prec_1, \prec_2}$ we put
$$
t_v = \cfrac{\<\omega_\delta, v\>}{\<\omega_\delta, v\> -
  \<\tau_\epsilon, v\>} = \cfrac{1}{1 - \cfrac{\<\tau_\epsilon,
    v\>}{\<\omega_\delta, v\>}}.
$$
If
$u, v\in \partial(I) \cap C_{\prec_1, \prec_2}$ then 
${\<\omega_\delta, u\>}$, ${\<\omega_\delta, v\>} > 0$
and 
\begin{align*}\tag{$*$}
  t_u &< t_v &\iff\\
  \cfrac{\<\tau_\epsilon, u\>}{\<\omega_\delta, u\>} &<
  \cfrac{\<\tau_\epsilon, v\>}{\<\omega_\delta, v\>}&\iff\\
  \<\tau_\epsilon, \<\omega_\delta, v\> u\> &<
  \<\tau_\epsilon, \<\omega_\delta, u\> v\> &\iff\\
  \<\omega_\delta, v\> u &\prec_2 \<\omega_\delta, u\> v
\end{align*}
To evaluate $\prec_2$ above we see that
\begin{align*}\tag{$**$}
  \<\tau_i, \<\omega_\delta, v\> u\> &< \<\tau_i, \<\omega_\delta, u\>
  v\>&\iff\\ 
  \<\omega_\delta, \<\tau_i, u\> v\> &< \<\omega_\delta, \<\tau_i, v\>
  u\>&\iff\\ 
  \<\tau_i, u\> v &\prec_1 \<\tau_i, v\> u
\end{align*}
for $i = 1, \dots, n$. Let $T$ denote the matrix whose rows are
$\tau_1, \ldots, \tau_n$. By choosing $\delta$
and $\epsilon$ generically as above it follows that
$$
t_u < t_v \iff
T u v^t \prec_1 T v u^t
$$
where $\prec_1$ above refers to the lexicographic extension of
$\prec_1$ on $\ZZ^n$ to $\ZZ^n \times \cdots \times \ZZ^n$. Here,
$Tuv^t$ and $Tvu^t$ are $n \times n$ matrices and we need to compare
their rows.  
Notice that the comparison between
$t_u$ and $t_v$ does not involve $\delta$ and $\epsilon$ but only the
term orders $\prec_1$ and $\prec_2$.  This leads us to define the {\it
facet preorder\/} $\prec$ by
\begin{equation}\label{facetpo}
u \prec v \iff t_u < t_v\iff T u v^t \prec_1 T v u^t 
\end{equation}
for $u, v\in \partial(I)\cap C_{\prec_1, \prec_2}$.
\begin{Example}\label{Example:facetpreorder}
  Continuing the setup in Example \ref{Example2}, if $u = (u_1, u_2)$
  and $v = (v_1, v_2)$, then
$$
T = 
\begin{pmatrix}
1 & 0\\
0 & 1
\end{pmatrix}
$$
and the facet preorder $\prec$ is given by
$$
u\prec v \iff (u_1 v \prec_1 v_1 u) \vee (
(u_1 v = v_1 u) \wedge (u_2 v \prec_1 v_2 u)).
$$
\end{Example}

If $t_u = t_v$ then $T u v^t = T v u^t$ and $u v^t = v u^t$ since
$T$ is an invertible matrix. The identity $u v^t = v u^t$ implies that
$u$ and $v$ are collinear. Since $u$ and $v$ lie in the same half space,
$u$ is a positive
multiple of $v$. 

This has
the nice consequence that the line $\omega(t)$ between $\omega_\delta$
and $\tau_\epsilon$ intersects the cones in the Gr\"obner fan in
dimension $\geq n-1$.  Consider the typical situation, where $v \in C
= C(v, v_1, \dots, v_m)$ is chosen to minimize $t_v$ as in the
Gr\"obner walk. Then $\omega(t_v)$ is on a proper face $F$ of $C^\vee$.
Since $t_v = t_u$ implies that $u$ is a positive multiple of $v$ for
$u\in \{v_1, \dots, v_m\}$, we conclude that $\dim F = n - 1$ i.e. $F$
is a facet.

The facet preorder $\prec$ defined in (\ref{facetpo}) may be inserted in 
the classical Gr\"obner walk algorithm giving the 
{\it generic Gr\"obner walk algorithm\/}
completely removing the numerical dependence on the line $\omega(t)$. 
Below, $-\infty (\infty)$ denotes a vector strictly 
smaller (larger) than the vectors in
$\partial(I)\cap C_{\prec_1, \prec_2}$.
 
\

\noindent
{\bf INPUT}: Marked reduced Gr\"obner basis $G$ for $I$ over a term
order $\prec_1$ and a term order $\prec_2$ (the facet preorder $\prec$ is 
given as in (\ref{facetpo}) using $\prec_1$ and $\prec_2$).

\

\noindent
{\bf OUTPUT}: Reduced Gr\"obner basis for $I$ over $\prec_2$.

\begin{enumerate}[(i)] 
\item
$w = -\infty$.
\item \label{Walk:Repeatg}
{\bf Compute\_last\_w}. 
If $w = \infty$ output $G$ and halt.
\item
Compute generators $\ini_\omega(G) = \{\ini_\omega(g)\mid g\in G\}$ for 
$\ini_\omega(I)$ as
$$
\ini_\omega(g) = a^u x^u + \sum_{v\in S_g} a_v x^v,
$$ 
where $S_g = \{v\in \supp(g)\setminus \{u\} \mid u - v \prec w, w\prec u - v\}$ and $a_u x^u$ is the marked term of $g\in G$.
\item Compute reduced Gr\"obner basis $H$ for $\ini_\omega(I)$
  over $\prec_2$ and mark $H$ according to $\prec_2$.
\item
Let
$$
H' = \{f - f^G \mid f\in H\}.
$$
Use marking of $H$ to mark $H'$.
\item
Autoreduce $H'$ and put $G = H'$.
\item
Repeat from (\ref{Walk:Repeatg}).
\end{enumerate}

\

\noindent
{\bf Compute\_last\_w:}

\begin{enumerate}
\item
Let 
$V := \{v\in \partial(G) \cap C_{\prec_1, \prec_2} \mid w \prec v\}$.
\item
If $V = \emptyset$, put $w = \infty$ and return. 
\item
Let $w := \min_\prec\{v | v\in V\}$ and return.
\end{enumerate}

\subsection{Variations on the generic Gr\"obner walk}
\

Several variations on the generic Gr\"obner walk are possible. In many
cases generators for an ideal are given which form a natural Gr\"obner
basis with respect to a specific weight vector. This happens for
example in implicitization problems with polynomials $y_1 - f_1,
\dots, y_m - f_m$, where $f_i$ are polynomials in $x_1, \dots, x_n$
for $i=1, \dots, m$. These polynomials form a Gr\"obner basis over a
vector $\omega$ assigning zero weights to $x_1, \dots, x_n$ and
positive weights to $y_1, \dots, y_m$. In this case one only needs to
work with $\omega$ and perturbations $\tau_\epsilon$ of the target
vector.  One may also truncate the facet preorder $\prec$ (to get a
{\it face preorder\/}) using only parts $(\omega_1, \dots, \omega_p)$
and $(\tau_1, \dots, \tau_q)$ of the $\QQ$-bases $\omega$ and $\tau$.
This leads to an analogue of the perturbation degree $(p,q)$-walk
defined in \cite{AGK}.

\section{An introductory example}\label{Section:Example}

We illustrate the generic Gr\"obner walk 
by a detailed example in the two 
dimensional case. For a given polynomial $f\in R$ we let 
$L_G(f) = f - f^G$, 
where $G$ is a marked Gr\"obner basis (markings are underlined).
Let $\prec_1$ denote degree (reverse) lexicographic
order and $\prec_2$ lexicographic order with $y\prec_{1,2} x$.
The facet preorder $\prec$ is given as in Example 
\ref{Example:facetpreorder}.
Consider the ideal 
$$
I = \<x^2 - y^3, x^3 - y^2 - x\>\subset \QQ[x, y].
$$
Initially we put  
$$
G=\{\underline{y^3} - x^2, \underline{x^3} - y^2 - x\},
$$
where the initial terms over $\prec_1$ are marked. 
Clearly $G$ is the reduced Gr\"obner 
basis for $I$ over $\prec_1$. The Gr\"obner cone is
given by
$$
C_{\prec_1}(I) = C(\{(-2, 3)\}\cup \{(3,-2)\})^\vee \cap \RR^2_{\geq 0}.
$$
In this case $(3, -2)\not\in C_{\prec_1, \prec_2}$ and 
$V = \{(-2, 3)\}$. So the first facet ideal is  
$\<\underline{y^3}-x^2, x^3\>$. The reduced Gr\"obner basis for this ideal over  
$\prec_2$ is $\{\underline{x^2} - y^3, xy^3, y^6\}$ and the lifting step 
is given 
by
\begin{align*}
L_G(x^2 - y^3) &= x^2 - y^3\\
L_G(xy^3) &= xy^3 - y^2 - x\\
L_G(y^6) &= y^6-xy^2-x^2.
\end{align*}
Our new marked reduced Gr\"obner basis is 
$$
G = \{\underline{x^2} - y^3, \underline{xy^3} - y^2 - x, 
\underline{y^6} - xy^2 - y^3\}.
$$
Since $w = (-2,3) \prec (-1,4)$ it follows that $V = \{(-1,4)\}$ and 
the next facet ideal is $\<x^2, xy^3, \underline{y^6} - xy^2\>$ 
with reduced Gr\"obner basis
$\{x^2, \underline{xy^2}-y^6, y^7\}$ over $\prec_2$. Since
\begin{align*}
L_G(x^2) &= x^2 - y^3\\
L_G(xy^2-y^6) &= xy^2 - y^6 + y^3\\
L_G(y^7) &= y^7 - y^4 - y^2 - x
\end{align*}
our new marked reduced Gr\"obner basis is 
$$
G = \{\underline{x^2} - y^3, \underline{xy^2} - y^6 + y^3, 
\underline{y^7} - y^4 - y^2 - x\}.
$$
Since $w = (-1,4)\prec (-1, 7)$ we get $V = \{(-1, 7)\}$ and the next facet 
ideal is $\<x^2, x y^2, \underline{y^7} - x\>$ with reduced 
Gr\"obner basis $(y^9, \underline{x}-y^7)$ over $\prec_2$. Here
\begin{align*}
L_G(y^9) &= y^9 - 2 y^6 - y^4 + y^3\\
L_G(x - y^7) &= x - y^7 + y^4 + y^2.
\end{align*}
The new marked reduced Gr\"obner basis is
$$
G = \{\underline{y^9} - 2 y^6 - y^4 + y^3, \underline{x} - y^7 + y^4 + y^2\}.
$$
Since $V = \emptyset$ in this case, the generic Gr\"obner walk halts and $G$ 
is the reduced Gr\"obner basis for $I$ over $\prec_2$.

\section{Computational experience for lattice ideals}
\label{Section:Computations}

In this section we report briefly on computations using 
the implementation {\tt GLATWALK} \cite{LGW} of the 
Buchberger algorithm and generic Gr\"obner walk for lattice ideals. Not surprisingly 
the walk performs best when initial and target vectors are close. 
An ideal situtation where this arises seems to come from a special case of
feasibility in integer programming. Consider natural numbers
$a_1, \dots, a_n\in \NN$. Given $b\in \NN$ decide if the equation
\begin{equation}\label{F1}
x_1 a_1 + \cdots + x_n a_n = b
\end{equation}
has a solution $x_1, \dots, x_n\in \NN$ and find
it if so. Adjoining the extra variable $t$ we seek to
minimize $t$ subject to
\begin{equation}\label{F2}
t + x_1 a_1 + \cdots + x_n a_n = b
\end{equation}
and $t, x_1, \dots, x_n \geq 0$. We denote this integer programming
problem $IP_{A, \tau}(b)$, where $A$ is the $1 \times (n+1)$ matrix
$(1\, a_1\, \dots \, a_n)$ and $\tau = (1, 0, \dots, 0)\in \NN^{n+1}$.
This problem has a trivial feasible
solution: $t=b,  x_1= \cdots = x_n=0$.
Now we may apply standard
algebraic techniques in integer programming (cf.~\cite{CT} and \cite{Rekha}) and form 
the toric ideal 
\begin{equation}\label{IA}
I_A = (x_1 - t^{a_1}, \dots, x_n - t^{a_n})\subset \QQ[t, x_1, \dots, x_n].
\end{equation}
A Gr\"obner basis $G_\tau$ for $I_A$ with respect to $\tau$ is a 
test set for the integer programming problems $IP_{A, \tau}(b)$, where 
$b$ varies and an optimal solution to (\ref{F2}) is the exponent
of the normal form of $t^b$ with respect $G_\tau$ thereby solving (\ref{F1}).

It is important to observe that the generating set for $I_A$ in (\ref{IA}) 
already is a Gr\"obner basis $G_\sigma$ for $I_A$ with respect to the vector
$\sigma = (-1, 0, \dots, 0)$. In the
following section we report on computational results in
computing $G_\tau$ using the generic walk to go from
$\sigma$ to $\tau$ compared with a direct computation with
Buchberger's algorithm. We use the programs {\bf walk} and
{\bf gbasis} of the program package {\tt GLATWALK}.

\subsection{Comparison with Buchberger's algorithm}

To walk from $\sigma$ to $\tau$ we break ties with the reverse
lexicographic order $<$ given by $t < x_1 < \cdots < x_n$ i.e.
we walk from the initial term order $<_\sigma$ to the target
term order $<_\tau$. The names of the computational examples in the 
following table refer to specific numbers $a_1, \dots, a_n$ as
in \S \ref{Section:Computations}. They can be found in
\cite{AL}. The timings below are in seconds and the computations were
carried out on a 1.6 GHz Pentium mobile with 1MB L2 Cache. 

\

\newcommand{\la}{\raggedleft\arraybackslash}

$$
\begin{tabular}{|>{\tt}p{3cm}|p{2.5cm}|p{2.5cm}|p{2.5cm}|p{2.5cm}|}
\hline
\hfil EXAMPLE\hfil & \hfil{\bf walk}\hfil & \hfil{\bf gbasis}\hfil & \hfil$|G_\sigma|$\hfil & \hfil$|G_\tau|$\hfil\\
\hline
cuww1 & \la 1.1 &  \la 17.7 & \la 5 & \la 7343\\
\hline
cuww2 & \la 11.4 & \la 2.4 & \la 6 &\la 2472\\
\hline
cuww3 & \la 24.4 & \la 9.5 & \la 6 & \la 4888\\
\hline
cuww4 & \la 1.2  & \la 21.3 & \la 7 & \la 7937\\
\hline
cuww5 & \la 7.9  & \la 1.3 & \la 8 & \la 1724\\
\hline
prob1 & \la 0.1  & \la 0.1 & \la 8 & \la 410\\
\hline
prob2 & \la 0.0  & \la 0.0 & \la 8 & \la 142\\
\hline
prob3 & \la 0.1  & \la 0.1 & \la 8 & \la 425\\
\hline
prob4 & \la 0.1  & \la 0.2 & \la 8 &  \la 757\\
\hline
prob5 & \la 0.2  & \la 0.1 & \la 8 & \la 516\\
\hline
prob6 & \la 0.1  & \la 0.5 & \la 10 & \la 1035\\
\hline
prob7 & \la 0.1  & \la 0.1 & \la 10 & \la 461\\
\hline
prob8 & \la 0.2 & \la 0.1 & \la 10 & \la 558 \\
\hline
prob9 & \la 0.0 & \la 0.0 & \la 10 & \la 270\\
\hline
prob10 & \la 0.6 & \la 2.5 & \la 10 & \la 2416\\
\hline
\end{tabular}
$$

\

\

In the problems {\tt cuww1, cuww4} and {\tt prob10} the initial and target vectors are
separated by less than $10$ Gr\"obner
cones in the Gr\"obner fan . This leads to surprisingly fast computation of 
relatively large 
Gr\"obner bases. It would be interesting to further explore the 
efficiency of the generic Gr\"obner walk in solving Frobenius
problems more general than the Aardal-Lenstra knapsack problems.

\section{Concluding Remarks}

The strength of the generic walk is that it is a completely
deterministic algorithm avoiding the inherent instability
of explicit numerical computation. 

The most recent version of the program {\tt 4ti2} (see \cite{H}) employs
a projection algorithm for computing Gr\"obner bases of
toric ideals before using the usual Buchberger algorithm.
This addition turns out to be a crucial optimization. It would
be interesting to use a similar projection algorithm before
using the walk. Our preliminary experiments indicate that
the walk is a little faster than the projection algorithm
in the cases where it is significantly faster than the
Buchberger algorithm. The current available implementations
of Buchberger's algorithm like CoCoa, Singular, Macaulay2 etc.~are
slowed down by a significant factor compared to
specialized integer vector implementations for lattice ideals.

While the generic Gr\"obner walk is presented here as a technique
to compute a Gr\"obner basis efficiently, one can use it for
walking the entire Gr\"obner fan systematically.  In fact,
a recent paper \cite{FJT} presents an algorithm based on both
the generic Gr\"obner walk and the reverse search technique
to list all Gr\"obner bases of a general polynomial ideal.
In short, it reverses the generic Gr\"obner walk in all possible
ways from the lexicographic basis to reach all other bases.
Obviously, such an exhaustive search requires an enormous amount
of computational effort, and the symbolic perturbation turns
out to be essential for this purpose.

\bibliographystyle{amsplain}

\end{document}